\documentclass[12pt, leqno, twoside]{amsart}
\usepackage{amsmath, amsthm, amscd, amsfonts, amssymb, graphicx, color}
\usepackage{mathrsfs}
\usepackage[bookmarksnumbered, colorlinks, plainpages]{hyperref}
\hypersetup{colorlinks=true,linkcolor=red, anchorcolor=green, citecolor=cyan, urlcolor=red, filecolor=magenta, pdftoolbar=true}
\usepackage[]{draftwatermark}
\SetWatermarkText{Draft}
\SetWatermarkScale{8}
\newcommand{\disp}{\displaystyle} 
\usepackage{enumerate}
\newtheorem{theorem}{Theorem}[section]
\newtheorem{lemma}[theorem]{Lemma}
\newtheorem{proposition}[theorem]{Proposition}
\newtheorem{corollary}[theorem]{Corollary}
\theoremstyle{definition}
\newtheorem{definition}[theorem]{Definition}

\theoremstyle{remark}
\newtheorem{remark}[theorem]{Remark}
\numberwithin{equation}{section}

\begin{document}

\setcounter{page}{1}

\title[Variable exponent Picone identity and $p(x)$ sub-Laplacian...]{Variable exponent Picone identity and $p(x)$ sub-Laplacian first eigenvalue  for general vector fields}

\author[A. Abolarinwa]{Abimbola Abolarinwa}

\author[A. Ali]{Akram Ali}

%\author[]{}

\address{Department of Mathematics,
University of Lagos, Akoka, Lagos State,
Nigeria, 
and}
\address{Department of Mathematics, Logic and Discrete Mathematics, 
Ghent University,  Ghent, Belgium.}
\email{\textcolor[rgb]{0.00,0.00,0.84}{a.abolarinwa1@gmail.com, aabolarinwa@unilag.edu.ng}}

\address{Department of Mathematics,
College of Science, King Khalid University, 9004 Abha, Saudi Arabia.}
\email{\textcolor[rgb]{0.00,0.00,0.84}{aakramali@ufam.edu.br; akramali133@gmail.com}}

% \date{June 20, 2022}

\begin{abstract}
In this paper,  we establish a new generalized nonlinear variable exponent Picone identities for $p(x)$-sub-Laplacian.  As applications we prove uniqueness, simplicity, monotonicity and isolatedness of the  first nontrivial Dirichlet eigenvalue of $p(x)$-sub-Laplacian with respect to the general vector fields.  Further applications yield Hardy type inequalities and Caccioppoli estimates with variable exponents.
\end{abstract}
\maketitle
 
 \tableofcontents
 
%%%%%%%%%%%%%%%%%%%%%%%%%%%%%%%%%%%%%%%%%%%%%%
%%%%%%%%%%%%%%%%%%%%%%%%%%%%%%%%%%%%%%%%%%%%%%
%%%%%%%%%%%%%%%%%%%%%%%%%%%%%%%%%%%%%%%%%%%%%% 
\section{Preliminaries}
\subsection{Introduction}
This paper is concerned with variable exponent Picone identity in the context of sub-Riemannian geometry. We derive a nonlinear Picone identity which allows us to study some qualitative properties of the principal eigenvalue of $p(x)$-sub-Laplacian with respect to the general vector fields on smooth manifolds. As by-products,  we also  derive Hardy type inequalities and Caccioppoli estimates with variable exponents. These results are appearing for the first time, even in the Euclidean setting. In recent years, several authors have devoted their researches towards the study of variable exponent elliptic equations and systems with $p(x)$-growth condition in Euclidean setting with many interesting results \cite{Alv,Deng,Fan,FZZ,HHV,MV}.  Models involving $p(x)$-growth condition arise from physical processes such as nonlinear elasticity theory,  electrorheological fluids,  image processing, etc \cite{AMS,AMSo,Ru}.  It has been observed that $p(x)$-Laplacian is similar in many respect to the classical $p$-Laplacian ($p$-constant) but it lacks certain vital properties such as homogeneity. This therefore makes the nonlinearity so much complicated and many of known approaches to $p$-Laplacian can no longer hold for $p(x)$-Laplacian. It is interesting to consider $p(x)$-Laplacian in the sub-elliptic setting and investigate which of the known results for $p$-constant hold for variable exponents.

\subsection{$p(x)$-Sub-Laplace operator and eigenvalues}
Let $M$ be an $n$-dimensional smooth manifold equipped with a volume form $dx$ and $\{X_k\}_{k=1}^N$, $n\ge N$, be a family of vector fields defined on $M$.   Consider the operator
$$\mathscr{L}_X := \sum_{k=1}^N X^*_kX_k,$$
which is a second-order differential operator usually called canonical sub-Laplacian.  This operator is related to the operator for the sum of squares of vector fields and it is well known to be locally hypoelliptic if the commutators of the vector fields $\{X_k\}_{k=1}^N$  generate the tangent space of $M$  as the Lie algebra, due to H\"ormander's pioneering work \cite{Hom}. We denote the horizontal gradients for general vector fields by 
$$\nabla_X = (X_1, \cdots, X_N) \ \ \ \text{and} \ \  \ \nabla_X^* = (X_1^*, \cdots, X_N^*),$$
where $X_k$ and its formal adjoint $X_k^*$ are respectively given by
$$X_k =  \sum_{j=1}^n a_{kj}(x)\frac{\partial}{\partial x_j} \ \ \text{and}\ \ X_k^* =  -\sum_{j=1}^n \frac{\partial}{\partial x_j} (a_{kj}(x)), \ \  k= 1, \cdots, N.$$
There are numbers of examples of sub-manifolds where vector fields can be defined.  For examples,  we list among others, the Carnot groups,  Heisenberg groups,  Engel groups,  and Grushin plane (which does not even posses a group structure).  Interested readers can see the book \cite{RS} for more examples and detail discussions on the sub-Laplacian  and its various extensions in each case.  In the case $M=\mathbb{R}^n$, then  $dx$ is the Lebesgue measure,  $\nabla_X=\nabla$ and $\mathscr{L}_X =\Delta$ are the usual Euclidean gradient and Laplacian, respectively.  

Let $p:\bar{\Omega} \to \mathbb{R}$ be a continuous function and $p(x)>1$ for $x\in \bar{\Omega}\subset M$. We  define  the $p(x)$-sub-Laplacian  for general vector fields on $M$ by the formula
$$\mathscr{L}_pu:= \nabla^*_X(|\nabla_X u|^{p(x)-2}\nabla_Xu),$$
where $u$ is a smooth function.  If $p(x)=p$ ($p$=constant),  the operator $\mathscr{L}_pu$ becomes the $p$-sub-Laplacian, $ \nabla^*_X(|\nabla_X u|^{p-2}\nabla_Xu)$ and $|x|$ stands for the Euclidean length of $x=(x_1,\cdots,x_n)$.  

As mentioned earlier,  various partial differential equations with variable exponent growth condition have appeared in literature (see \cite{Alv,Deng,Fan, FZZ,HHV,MV} for instance), but there is scarcity of such mathematical models in the subelliptic setting.  In this paper however we shall consider the indefinite weighted  Dirichlet eigenvalue problem for $p(x)$-sub-Laplacian on $\Omega\subset M$, $p(x)>1$,
\begin{align}\label{eq1}
\left.
\begin{array}{ll}
-\nabla^*_X(|\nabla_X u|^{p(x)-2}\nabla_Xu) = \lambda g(x)|u|^{p(x)-2} u,  & \ x \in \Omega,\\
\ \   \ u>0,  &\ x \in \Omega,\\
\ \ \ u=0,   &\ x \in \partial\Omega,
\end{array}
\right.
\end{align}
and discuss some properties of the eigenvalue $\lambda \in \mathbb{R}^+$ and the corresponding eigenfunction $u(x)$ in certain Sobolev spaces with variable exponents \cite{DHHR,CFRW,FZ}. It is well known in the classical setting ($p(x)=p$-constant and $M=\mathbb{R}^n$) that Problem \eqref{eq1} possesses a closed set of nondecreasing sequence of nonnegative eigenvalues $\{\lambda_k\}$ which grows to $+\infty$ as $k\to +\infty$, and that the first nonzero eigevalue is simple and isolated.  Due to some complication in the nonlinearities in $p(x)$-Laplacian and inhomogeneity  of the  corresponding variable exponent norm, some of the results in the classical case may not hold or rather under restrictive assumptions.  In \cite{FZZ}, the authors studied \eqref{eq1} (with $g(x)=1$,  $M=\mathbb{R}^n$) and showed the existence of infinitely many eigenvalues and established some sufficient condition for the infimum of the spectrum (called the principal eigenvalue), 
$$\lambda_{1,p} = \inf _{u\ne 0}\frac{\int_\Omega |\nabla u|^{p(x)}dx}{ \int_\Omega |u|^{p(x)}dx}, \ \ \ p(x)>1,$$
 to be zero and positive, respectively.  The properties that $\lambda_{1,p}>0$ is very useful in analysis and applications.  Motivated by \cite{FZZ}, we are able to assume the existence of $\lambda_{1,p}>0$ for \eqref{eq1} and proved its uniqueness, monotonicity, simplicity and isolatedness. The  variable exponent Picone identity (discussed in Section \ref{sec2}) plays a crucial role in our proofs. 

\subsection{Picone identities} Picone identity is a very useful tool in the study of qualitative properties of solutions of differential equations,  and for this,  several  linear and nonlinear Picone type identities have been derived to handle differential equations of various type.  Picone identity was originally developed by Mauro Picone in 1910 to prove Sturm Comparison principle and oscillation theory for a system of differential equations. This identity was later extended to  partial differential equation involving Laplacian by Allegretto  \cite{Al1} and  $p$-Laplacian by Allegretto and Huang \cite{AH} to establish among others,  existence and nonexistence of positive solutions,  Sturmian comparison principle,  Liouville type theorems, Hardy inequalities and some profound results involving $p$-Laplace equations and systems. Precisely, Allegretto  \cite{Al1} proved that, for nonnegative differentiable functions $u$ and $v$ with $v \neq 0$, the following formula
\begin{align}\label{e11}
 |\nabla u|^2+\frac{u^2}{v^2}|\nabla v|^2-2\frac{u}{v}\nabla u\nabla v= |\nabla u|^2-\nabla\left(\frac{u^2}{v}\right)\nabla v \ge 0
\end{align}
holds.  Allegretto and Huang \cite{AH} extended  \eqref{e11} to handle $p$-Laplace equations and eigenvalue problems involving $p$-Laplacian. Their identity reads as follows, for $u\ge 0$, $v>0$, then 
\begin{align}\label{e12}
|\nabla u|^p+(p-1)\frac{u^p}{v^p}|\nabla v|^p &-p\frac{u^{p-1}}{v^{p-1}}|\nabla v|^{p-2}\nabla v\nabla u = R_p(u,v),
\end{align}
where
\begin{align*}
R_p(u,v):=|\nabla u|^p-\nabla\left(\frac{u^p}{v^{p-1}}\right)|\nabla v|^{p-2}\nabla v \ge 0.
\end{align*} 
Several extensions and generalization of Picone identity have been established in order to handle more general elliptic operators.  Tyagi \cite{Ty} and Bal \cite{Ba} established nonlinear versions of \eqref{e11} and its $p$-Laplace analogue \eqref{e12}, respectively, with several applications, (see also \cite{DT,Fe,Tir}).  For other interesting extension of Picone type identities one can find \cite{Ja3,Ja4} (for Finsler $p$-Laplacian  with application to Caccioppoli inequality), \cite{RSS1,RSS2,RSS3} (for general vector fields and $p$-sub-Laplacian with applications to Grushin plane, Heisenberg group, Stratified Lie groups), \cite{NZW} (for $p$-sub-Laplacian on Heisenberg group and applications to Hardy inequalities), \cite{A1,SY} (for nonlinear Picone identities for anisotropic $p$-sub-Laplacian and $p$-biLaplacian with applications to horizontal Hardy inequalities and weighted eigenvalue problem on Stratified Lie groups).

Allegretto \cite{Al3} established variable exponent Picone type identity for differentiable functions $v>0$,  $0\le u\in C^\infty_0(\Omega)$, $\Omega\subset\mathbb{R}^n$ with $n\ge 1$ and continuous $p(x)>1$ as follows:
\begin{align}\label{e14}
\frac{|\nabla u|^{p(x)}}{p(x)} &- \nabla\left[\frac{u^{p(x)}}{p(x)v^{p(x)-1}} \right]|\nabla v|^{p(x)-2}\nabla v \nonumber\\
& = \frac{|\nabla u|^{p(x)}}{p(x)} - \left(\frac{u}{v}\right)^{p(x)-1}|\nabla v|^{p(x)-2}\nabla v\nabla u +\frac{p(x)-1}{p(x)} \left(\frac{u}{v}|\nabla v| \right)^{p(x)}\\
&  \hspace{1cm} + \frac{1}{p(x)} \frac{u^{p(x)}}{v^{p(x)-1}}|\nabla v|^{p(x)-2} \left[\frac{1}{p(x)} - \ln\left(\frac{u}{v}\right)\right]\nabla v\nabla p(x) \ge 0\nonumber 
\end{align}
on the assumption that $\nabla v\nabla p(x)=0$.  He used the inequality to prove Barta theorem and some other results.  Later, Yoshida \cite{Yo1} (see also \cite{Yo2,Yo3}) established  similar Picone identities for quasilinear  and half-linear elliptic equations involving $p(x)$-Laplacian and pseudo $p(x)$-Laplacian, and consequently developed Sturmian comparison theory. Most recently, Feng and Han \cite{FH}, motivated by Allegretto \cite{Al3} proved a modified form of \eqref{e14} and showed that 
\begin{align}
|\nabla u|^{p(x)} - \nabla\left(\frac{u^{p(x)}}{v^{p(x)-1}} \right)|\nabla v|^{p(x)-2}\nabla v \ge 0
\end{align}
if $\nabla v\nabla p(x)=0$ a.e in $\Omega$, with equality if and only if $\nabla(u/v)=0$ in $\Omega$. They proved monotonicity of principal eigenvalue $\lambda_{1,p}$ and a variable  exponent Barta inequality for  $p(x)$-Laplacian in the form
\begin{align*}
\lambda_{1,p} \ge \inf_{x\in \Omega}\left[\frac{\Delta_pv}{v^{p(x)-1}}\right],  \ \ \ \Omega\subset \mathbb{R}^n,
\end{align*}
where $\Delta_p:=-\nabla(|\nabla v|^{p(x)-2}\nabla v)$, on the assumption that $\nabla v\nabla p(x)=0$.

\subsection{Variable  exponent functional spaces}
In order to discuss generalized solutions, we need some concepts from the theory of variable Lebesgue and Sobolev spaces. Detailed description of these spaces can be found in \cite{CFRW,DHHR,FZ}. 

Let $\Omega\subset M$ be an open domain and $E(\Omega)$ denotes the set of all equivalence classes of measurable real-valued functions defined on $\Omega$ being equal almost everywhere. 

\begin{definition}
The variable exponent Lebesgue space $L^{p(\cdot)}(\Omega)$ is defined as 
$$L^{p(\cdot)}(\Omega) = \left\{u\in E(\Omega) : \int_\Omega|u(x)|^{p(\cdot)} dx < \infty\right\}$$
equipped with the (Luxemburg) norm
$$\|u\|_{L^{p(\cdot)}(\Omega)}= \inf\left\{t>0: \int_\Omega \left| \frac{u(x)}{t}\right|^{p(x)} dx \le 1 \right\}.$$
\end{definition}

Consider the functional (also called the $\rho$-modular) on $L^{p(\cdot)}(\Omega)$, which is the mapping $\rho_{p(\cdot)}(u): L^{p(\cdot)}(\Omega)\to \mathbb{R}$, and defined by 
$$L^{p(\cdot)}(\Omega):=\int_\Omega|u(x)|^{p(x)} dx.$$
The following proposition contains vital results in the study of variable Lebesgue space.   We suppose a continuous function $p:\bar{\Omega}\to \mathbb{R}^+$, $p(x)>1$ is such that
$$1<p^-:=ess\inf_{x\in\bar{\Omega}} p(x)\le p(x) \le p^+:= ess\sup_{x\in\bar{\Omega}} p(x)<\infty.$$
\begin{proposition}\cite{CFRW,DHHR,FZ}\\
Denote $\|u\|_{p(x)} :=\|u\|_{L^{p(x)}(\Omega)}$. For any $u, u_m \in L^{p(x)}(\Omega)$, where $m=1,2,\cdots,$ the following statements are true:
\begin{enumerate}[1.]
\item $\|u\|_{p(x)}<1 (=1 \ \text{or}\ >1)$ if and only if $\rho_{p(x)}(u) <1 (=1 \ \text{or}\ >1)$;
\item If $\|u\|_{p(x)} \le 1$ then $\|u\|_{p(x)}^{p^+} \le \rho_{p(x)}(u)\le \|u\|_{p(x)}^{p^-}$;
\item If $\|u\|_{p(x)} > 1$ then $\|u\|_{p(x)}^{p^-} \le \rho_{p(x)}(u)\le \|u\|_{p(x)}^{p^+}$;
\item $\|u_m-u\|_{p(x)} \to 0$ if and only if $\rho_{p(x)}(u_m-u) \to 0$;
\item $\min\{ \|u\|_{p(x)}^{p^-}, \|u\|_{p(x)}^{p^+}\} \le \rho_{p(x)}(u) \le \max\{ \|u\|_{p(x)}^{p^-}, \|u\|_{p(x)}^{p^+}\} $.
\end{enumerate}
\end{proposition}

\noindent The following generalized H\"older's inequality can be used to define equivalent norms. 

\begin{proposition} (H\"older's inequality \cite{CFRW,DHHR})\\
Let $\frac{1}{p(x)}+\frac{1}{p'(x)}=1$ a.e. on $\Omega$, then for all $u\in L^{p(x)}(\Omega)$ and $v\in L^{p'(x)}(\Omega)$ we have $uv\in L^1(\Omega)$ and 
$$\int_\Omega|u(x)v(x)|dx \le \left(1+\frac{1}{p^-}-\frac{1}{p^+} \right) \|u\|_{p(x)}\|v\|_{p'(x)}.$$
\end{proposition}

\begin{definition}
The variable exponent Sobolev space $W^{1, p(\cdot)}(\Omega)$ is defined as 
$$W^{1, p(\cdot)}(\Omega) = \{u\in L^{p(\cdot)}(\Omega):  |\nabla_Xu| \in L^{p(\cdot)}(\Omega)\}$$
equipped with the norm
$$\|u\|_{W^{1, p(\cdot)}(\Omega)}= \|u\|_{L^{p(\cdot)}(\Omega)} + \|\nabla_Xu\|_{L^{p(\cdot)}(\Omega)}.$$

The space $W_0^{1, p(\cdot)}(\Omega)$ is defined as the closure of $C^\infty_0(\Omega)$ in $W^{1, p(\cdot)}(\Omega)$ with respect to the norm
$$\|u\|_{W_0^{1, p(\cdot)}(\Omega)}=  \|\nabla_Xu\|_{L^{p(\cdot)}(\Omega)}.$$
\end{definition}
\noindent It can be easily proved that $L^{p(\cdot)}(\Omega)$,  $W^{1, p(\cdot)}(\Omega)$ and $W_0^{1, p(\cdot)}(\Omega)$ are all separable and reflexive Banach spaces in their respectful norms if $1<\inf p(x)<\sup p(x)<\infty$ on $\Omega$.

\subsection{Plan of the paper}
In this paper, we derive new generalized variable exponent Picone type identities for general vector fields in the sub-Riemannian settings. The derived generalized identity contains some known Picone type identities in various settings as will be discussed in Section \ref{sec2}.  Consequently, we give several applications to qualitative properties of  the principal eigenvalue of $p(x)$-sub-Laplacian. Here, we are concerned with uniqueness, simplicity, monotonicity and isolatedenss of the  Dirichlet principal eigenvalue. These are discussed in Section \ref{sec3}.  Lastly, motivated by \cite{Ja4,RSS1}, we derive as a consequent of Picone identity,  sub-elliptic variable exponents Caccioppoli estimates in the form
\begin{align*}
\int_\Omega \phi^{p(x)}|\nabla_Xv|^{p(x)}dx  \le (p^+)^{p^+} \int_\Omega v^{p(x)} |\nabla_X\phi|^{p(x)} dx
\end{align*}
for every nonnegative test function $\phi \in C^\infty_0(\Omega)$, where $v$ is a sub-solution  in $\Omega\subset M$ and $p^+:=ess\sup p(x)$.  On the other hand, for $v$ positive $p(x)$-superharmonic functions,  we obtain a new version of  logarithmic Caccioppolli inequality 
\begin{align*}
\int_\Omega |\phi\nabla_X\log v|^pdx \le \left(\frac{p^+}{p^--1}\right)^{p^+}\int_\Omega  |\nabla_X\phi|^pdx.
\end{align*}

\section{Nonlinear variable exponent Picone identity} \label{sec2}
Here we give the statement and the proof of the nonlinear Picone identity with variable exponent, which is the main result of this section. First, we state some hypotheses as adopted in this section (and ofcourse throughout the paper) and Young's inequality in the forms that will be applied here and later.

Let $M$ be an $n$-dimensional smooth manifold and $\Omega$  any domain in $M$,  $p(x)>1$ is a continuous function on $\bar{\Omega}$, $p'(x) = 1/(p(x)-1)$ is H\"older conjugate to $p(x)$.
\begin{lemma}\label{lemY}
(Classical Young's inequality)
Let $s\ge 0$, $t\ge0$,  and $p(x)>1$ such that $1/p(x)+1/p'(x)=1$. There holds the inequality
\begin{align}\label{Y1}
st \le \frac{s^{p(x)}}{p(x)} + \frac{t^{p'(x)}}{p'(x)}
\end{align}
with equality if and only if $s^{p(x)}=t^{p'(x)}$.
\end{lemma}
\noindent Inequality \eqref{Y1} is the classical Young's inequality which can be varied in the following form.
\begin{lemma}\label{Y2}
(Modified Young's inequality)
Let $\Phi(x),\Psi(x)\ge 0$,   $p(x)>1$ such that $1/p(x)+1/p'(x)=1$ and $\varepsilon:\Omega\to \mathbb{R}^+$ be a continuous and bounded function.  There holds the inequality
\begin{align}\label{Y2}
\Phi\Psi^{p(x)-1} \le \frac{\Phi^{p(x)}}{p(x)\varepsilon(x)^{p(x)-1}} + \frac{p(x)-1}{p(x)}\varepsilon(x) \Psi^{p(x)}
\end{align}
for a.e. $x\in \Omega$ Furthermore, there is equality in \eqref{Y2} if and only if $\Phi =\varepsilon(x)\Psi$.
\end{lemma}
\proof Applying the classical Young's inequality \eqref{Y1} with
$$ s = \frac{\Phi}{\varepsilon(x)^{\frac{p(x)}{p(x)-1}}} \ \ \ \text{and}\ \ \ t = \left(\Psi \varepsilon(x)^{\frac{1}{p(x)}}\right)^{p(x)-1},$$
we have 
\begin{align*}
\Phi\Psi^{p(x)-1}& = \left(\frac{\Phi}{\varepsilon(x)^{\frac{p(x)}{p(x)-1}}}\right)  \left(\Psi \varepsilon(x)^{\frac{1}{p(x)}}\right)^{p(x)-1}\\
&\le   \frac{\Phi^{p(x)}}{p(x)\varepsilon(x)^{p(x)-1}} + \frac{p(x)-1}{p(x)}\left(\Psi \varepsilon(x)^{\frac{1}{p(x)}}\right)^{p(x)}.
\end{align*}
\qed

The next is the variable exponent Picone identity which is the main theorem in this section.
\begin{theorem}\label{Pic-thm}
 Let  $u\ge 0$ and $v>0$ be nonconstant differentiable functions a.e.  in $\Omega$. Suppose  $p:\bar{\Omega}\to (0,\infty)$ is a $C^1$-function for $p(x)>1$, and $f:(0,\infty)\to (0,\infty)$ is a $C^1$-function satisfying $f(y)>0$ and $f'(y)\ge (p(x)-1)\left[f(y)^{\frac{p(x)-2}{p(x)-1}}\right]$ for $y>0$. Define 
\begin{align}\label{e23}
L(u,v) &= |\nabla_Xu|^{p(x)}- \frac{u^{p(x)}\ln u}{f(v)} |\nabla_Xv|^{p(x)-2}\nabla_Xv\nabla_Xp(x) \nonumber \\
&-p(x)\frac{u^{p(x)-1}}{f(v)}  |\nabla_Xv|^{p(x)-2}\nabla_Xv \nabla_Xu 
+ \frac{u^{p(x)}f'(v)}{(f(v))^2}  |\nabla_Xv|^{p(x)}
\end{align}
and
\begin{align}\label{e24}
R(u,v) = |\nabla_Xu|^{p(x)}- \nabla_X\left(\frac{u^{p(x)}}{f(v)} \right)|\nabla_Xv|^{p(x)-2}\nabla_Xv.
\end{align}
Then
\begin{enumerate}
\item $L(u,v)=R(u,v)$.
\item Moreover $L(u,v)\ge 0$ if $\nabla_Xv\nabla_Xp(x)\equiv 0$.
\item Furthermore, $L(u,v)=0$ a.e. in $\Omega$ if and only if $\nabla_X(u/v)=0$ a.e. in $\Omega$.
\end{enumerate}
\end{theorem}

\proof
By direct computation we have 
\begin{align*}
R(u,v) & = |\nabla_Xu|^{p(x)}- \left(\frac{\nabla_X(u^{p(x)})}{f(v)} - \frac{u^{p(x)}\nabla_X(f(v))}{(f(v))^2} \right)|\nabla_Xv|^{p(x)-2}\nabla_Xv\\
&=  |\nabla_Xu|^{p(x)} - \frac{u^{p(x)}\ln u \nabla_Xp(x) +p(x)u^{p(x)-1}\nabla_Xu}{f(v)} |\nabla_Xv|^{p(x)-2}\nabla_Xv\\
& \hspace{2cm} + \frac{u^{p(x)}f'(v)}{(f(v))^2}  |\nabla_Xv|^{p(x)}\\
& = L(u,v),
\end{align*}
which proves $(1)$ of the theorem.

Next we verify $L(u,v)\ge 0$. Rewriting the expression for $L(u,v)$ as follows
\begin{align*}
L(u,v) &= |\nabla_Xu|^{p(x)}-p(x)\frac{u^{p(x)-1}}{f(v)}  |\nabla_Xv|^{p(x)-1}|\nabla_Xu | + \frac{u^{p(x)}f'(v)}{(f(v))^2}  |\nabla_Xv|^{p(x)} \\
& \ \ \ +p(x)\frac{u^{p(x)-1}}{f(v)} \left(|\nabla_Xv| |\nabla_Xu| - \nabla_Xv\nabla u \right) - \frac{u^{p(x)}\ln u}{f(v)} |\nabla_Xv|^{p(x)-2}\nabla_Xv\nabla_Xp(x)\\
& = p(x)\left(\frac{|\nabla_Xu|^{p(x)}}{p(x)} +\frac{p(x)-1}{p(x)} \left[\frac{(u|\nabla_Xv|)^{p(x)-1}}{f(v)} \right]^{\frac{p(x)}{p(x)-1}}\right) + \frac{u^{p(x)}f'(v)}{(f(v))^2}  |\nabla_Xv|^{p(x)}\\
& \ \ \ -(p(x)-1) \left[\frac{(u|\nabla_Xv|)^{p(x)-1}}{f(v)} \right]^{\frac{p(x)}{p(x)-1}} -p(x)\frac{u^{p(x)-1}}{f(v)}  |\nabla_Xv|^{p(x)-1}|\nabla_Xu | \\
&\ \ \ +p(x)\frac{u^{p(x)-1}}{f(v)} \left(|\nabla_Xv| |\nabla_Xu| - \nabla_Xv\nabla u \right)
- \frac{u^{p(x)}\ln u}{f(v)} |\nabla_Xv|^{p(x)-2}\nabla_Xv\nabla_Xp(x)\\
& = L_1(u,v) \ \ + \ \ L_2(u,v) \ \ + \ \ L_3(u,v)  \ \ + \ \ L_4(u,v),
\end{align*}
where
\begin{align*}
 L_1(u,v)&:= p(x)\left(\frac{|\nabla_Xu|^{p(x)}}{p(x)} +\frac{p(x)-1}{p(x)} \left[\frac{(u|\nabla_Xv|)^{p(x)-1}}{f(v)} \right]^{\frac{p(x)}{p(x)-1}}\right)\\
& \hspace{3cm} -p(x)\frac{u^{p(x)-1}}{f(v)}  |\nabla_Xv|^{p(x)-1}|\nabla_Xu |,
\end{align*}
\begin{align*}
 L_2(u,v)&:= \frac{u^{p(x)}f'(v)}{(f(v))^2}  |\nabla_Xv|^{p(x)} - (p(x)-1 )\left[\frac{(u|\nabla_Xv|)^{p(x)-1}}{f(v)} \right]^{\frac{p(x)}{p(x)-1}},
\end{align*}
\begin{align*}
 L_3(u,v)&:=  p(x)\frac{u^{p(x)-1}}{f(v)} \left(|\nabla_Xv| |\nabla_Xu| - \nabla_Xv\nabla u \right),
\end{align*}
\begin{align*}
 L_4(u,v)&:=  - \frac{u^{p(x)}\ln u}{f(v)} |\nabla_Xv|^{p(x)-2}\nabla_Xv\nabla_Xp(x).
\end{align*}
Applying the Young's inequality  \eqref{Y1}, choosing  $s=|\nabla_X u|$ and
$\disp t= \frac{(u|\nabla_Xv|)^{p(x)-1}}{f(v)}$, we obtain
\begin{align*}
p(x)\frac{u^{p(x)-1}}{f(v)} & |\nabla_Xv|^{p(x)-1}|\nabla_Xu | \\
& \le p(x)\left(\frac{|\nabla_Xu|^{p(x)}}{p(x)} +\frac{p(x)-1}{p(x)} \left[\frac{(u|\nabla_Xv|)^{p(x)-1}}{f(v)} \right]^{\frac{p(x)}{p(x)-1}}\right),
\end{align*}
implying that $L_1(u,v)\ge 0$ with equality if and only if there is equality in the Young's inequality, that is, $s=t^{\frac{1}{p(x)-1}}$.

Applying the assumption $f'(y)\ge (p(x)-1)\left[f(y)^{\frac{p(x)-2}{p(x)-1}}\right]$,  we have 
\begin{align*}
 \frac{u^{p(x)}f'(v)}{(f(v))^2}  |\nabla_Xv|^{p(x)} \ge  (p(x)-1 )\left[\frac{(u|\nabla_Xv|)^{p(x)-1}}{f(v)} \right]^{\frac{p(x)}{p(x)-1}},
\end{align*}
which implies that $L_2(u,v)\ge0$ with equality if and only if\\ $f'(y)= (p(x)-1)\left[f(y)^{\frac{p(x)-2}{p(x)-1}}\right]$. Clearly, $L_3(u,v)\ge 0$ by reverting to the inequality $|\nabla_Xv| |\nabla_Xu| - \nabla_Xv\nabla_X u\ge 0$. By the virtue of the assumption that $\nabla_Xv\nabla_Xp(x)\equiv 0$, we have also $L_4(u,v)\equiv 0$. Putting all of these together we obtain that $L(u,v)\ge 0$ a.e. in $\Omega$.

Observe that $L(u,v)=0$ holds if and only if 
\begin{align}\label{a23}
|\nabla_Xu| = \frac{u}{f(v)^{\frac{1}{p(x)-1}}} |\nabla_Xv|,
\end{align}
\begin{align}\label{a24}
f'(y)= (p(x)-1)\left[f(y)^{\frac{p(x)-2}{p(x)-1}}\right],
\end{align}
and
\begin{align}\label{a25}
|\nabla_Xv| |\nabla_Xu| =  \nabla_Xv\nabla_X u.
\end{align}
Upon solving for \eqref{a24} we get $f(v)=v^{p(x)-1}$. If $\nabla_X(u/v)=0$ then there exists a positive constant, say $\alpha>0$ such that $u=\alpha v$, then equality \eqref{a25} holds. Combining $f(v)=v^{p(x)-1}$ and $u=\alpha v$, then \eqref{a23} holds. 
We can now conclude that $L(u,v)=0$ implies $\nabla_X(u/v)=0$. Indeed, if $L(u,v)(x_0)=0$, $x_0\in \Omega$, there are two cases to consider, namely; the case $u(x_0)\neq 0$ and the case $u(x_0)=0$.\\
(a)\ If $u(x_0)\neq 0$, then $L(u,v)=0$ for all $x_0\in \Omega$, that is, $L_1(u,v)=0$, $L_2(u,v)=0$ and $L_3(u,v)=0$, and we conclude that \eqref{a23}, \eqref{a24} and \eqref{a25} hold, which when combined gives $u=\alpha v$ a.e. for some constant $\alpha>0$ and $\nabla_X(u/v)=0$ for all $x_0\in \Omega$.\\
(b)\ If $u(x_0) = 0$, we denote $\Omega^*=\{x\in \Omega: u(x)=0\}$, and suppose $\Omega^*\neq \Omega$. Here $u(x_0)=\alpha v(x_0)$ implies $\alpha=0$ since $u(x_0)=0$ and $v(x_0)>0$. By the first case (Case (a)) we know that $u(x)=\alpha v(x)$ and $u(x)\neq 0$ for all $x\in \Omega\setminus\Omega^*$, then it is impossible that $\alpha=0$. This contradiction implies that $\Omega^*=\Omega$.

\qed

\begin{remark}
Theorem \ref{Pic-thm} generalizes many known results. For examples:
 \begin{enumerate} 
 \item If $M=\mathbb{R}^n$ and $f(v)=v^{p(x)-1}$ in \eqref{e23} and \eqref{e24}. Then, we obtain the variable exponent Picone identity of Allegretto \cite{Al3} and Feng and Han \cite{FH}.
 \item If $p(x)=p$,  $f(v)=v^{p-1}$ in \eqref{e23} and \eqref{e24}, then our result covers Allegretto and Huang's \cite{AH} ($M=\mathbb{R}^n$), Niu, Zhang and Wang \cite{NZW} (Heisenberg group), Ruzhansky, Sabitbek and Suragan \cite{RSS1} (for general vector fields).
 \item If we allow $p(x)=p$ in \eqref{e23} and \eqref{e24},   we then recover Bal \cite{Ba} in the Euclidean setting and  Suragan and Yessirkegenov \cite{SY} in the setting of stratified Lie groups.
 \end{enumerate}
\end{remark}

\section{Applications}\label{sec3}

\subsection*{Eigenvalue problem for $p(x)$-sub-Laplacian}
Let $\Omega\subset M$ be a bounded domain with smooth boundary $\partial\Omega$. We suppose a continuous function $p:\bar{\Omega}\to \mathbb{R}^+$, $p(x)>1$ is such that
$$1<p^-:=ess\inf_{x\in\bar{\Omega}} p(x)\le p(x) \le p^+:= ess\sup_{x\in\bar{\Omega}} p(x)<\infty.$$
Now consider the indefinite weighted Dirichlet eigenvalue problem for $p(x)$-Laplacian
\begin{align}\label{e31}
\left.
\begin{array}{ll}
-\nabla^*_X(|\nabla_X u|^{p(x)-2}\nabla_Xu) = \lambda g(x)|u|^{p(x)-2} u,  & \ x \in \Omega,\\
\ \   \ u>0,  &\ x \in \Omega,\\
\ \ \ u=0,   &\ x \in \partial\Omega,
\end{array}
\right.
\end{align}
where $\Omega$ is as defined above,  $g(x)$ is a positive bounded function and $p:\bar{\Omega}\to (1,\infty)$ is a continuous function for $x\in \bar{\Omega}$.
\begin{definition}
Let $\lambda\in \mathbb{R}^+$ and $u\in W_0^{1,p(x)}(\Omega)$, the pair $(u,\lambda)$ is called a solution of \eqref{e31} if 
\begin{align}\label{e32}
\int_\Omega|\nabla_Xu|^{p(x)-2}\langle\nabla_Xu,\nabla_X\phi\rangle dx - \lambda \int_\Omega g(x)|u|^{p(x)-2} u\phi dx =0
\end{align}
for all $\phi \in W_0^{1,p(x)}(\Omega)$.
If $(u,\lambda)$ is a solution of \eqref{e31}, we call $\lambda$ an eigenvalue, and $u$ an eigenfunction corresponding to $\lambda$.

Similarly, by the sup-solution and sub-solution of \eqref{e31}, we mean the pair $(u,\lambda)$ such that 
\begin{align}\label{e33}
\int_\Omega|\nabla_Xu|^{p(x)-2}\langle\nabla_Xu,\nabla_X\phi\rangle dx - \lambda \int_\Omega g(x)|u|^{p(x)-2} u\phi dx \ge 0
\end{align}
and
\begin{align}\label{e34}
\int_\Omega|\nabla_Xu|^{p(x)-2}\langle\nabla_Xu,\nabla_X\phi\rangle dx - \lambda \int_\Omega g(x)|u|^{p(x)-2} u\phi dx \le 0
\end{align}
for all $\phi \in W_0^{1,p(x)}(\Omega)$, respectively.
\end{definition}
Denote the principal eigenvalue of \eqref{e31} (the least positive eigenvalue) by $\lambda_{1,p}:= \lambda_{1,p}(\Omega)$, clearly for the solution $(u,\lambda)$ and $u\neq 0$, we get
$$\lambda_{1,p} =\inf_{u\in \in W_0^{1,p(x)}(\Omega)\setminus\{0\}} \frac{\int_\Omega|\nabla_Xu|^{p(x)}dx}{\int_\Omega g(x)|u|^{p(x)}dx}.$$
In the case $p(x)=p$(constant), it is well known that $\lambda_{1,p}(\Omega)$ given above is the first eigenvalue of $p$-Laplacian (with $g(x)=1$,  $\Omega\subset\mathbb{R}^n$), which must be positive. But this is not true for general $p(x)$ in the sense that $\lambda_{1,p}$ may be zero \cite{FZ}. Nevertheless,  Fan, Zhang and Zhao in \cite{FZZ} have proved  the existence of infinitely many eigenvalues $p(x)$-Laplacian and established sufficient conditions for $\lambda_{1,p}(\Omega)>0$ (see also Franzina and  Lindqvist \cite{FL}).  Motivated by \cite{FZZ}, we are able to assume the existence of $\lambda_{1,p}>0$ in the rest of this section.   

In the rest of this section we are concerned with the indefinite weighted Dirichlet eigenvalue problem \eqref{e31} and discuss some properties of its solutions $v$ satisfying $\nabla_Xv\nabla_Xp(x)\equiv 0$ by the application of Picone identity in Theorem \ref{Pic-thm}. We remark that the results of this paper are classical in the sense that they have been established using different methods such variational approach (see \cite{Alv,FZZ,MV} for instance) where the condition $\nabla_Xv\nabla_Xp(x)\equiv 0$ is not required.

\subsection{Variable exponent Hardy type inequality}
\begin{proposition}\label{Pro33}
Let $\Omega\subset M$ be an open bounded domain. Suppose that a function $v\in C^\infty_0(\Omega)$ satisfies $\nabla_Xv\nabla_Xp(x)\equiv 0$ and 
\begin{align}\label{e35}
\left.
\begin{array}{ll}
-\mathscr{L}_pv= \mu a(x)f(v) & \ \  \text{in}\  \Omega,\\
\ \   \ v>0 &\ \  \text{in}\  \Omega,\\
\ \ \ v=0  &\ \  \text{on}\  \partial\Omega,
\end{array}
\right.
\end{align}
where $f:\mathbb{R}^+\to\mathbb{R}^+$ is $C^1$ and satisfies $f'(y)\ge(p(x)-1)\left[f(y)^{\frac{p(x)-2}{p(x)-1}} \right]$, $\mu>0$ is a constant, $a(x)$ is a positive continuous function. Then there holds
\begin{align*}
\int_\Omega|\nabla_Xu|^{p(x)}dx\ge \mu \int_\Omega a(x)|u|^{p(x)}dx
\end{align*}
for any $0\le u\in C^1_0(\Omega)$.
\end{proposition}

\proof
Since $v>0$ and solves \eqref{e35} in $\Omega$, that is, $v\in W_0^{1,p(x)}(\Omega)$. For a given a $\epsilon>0$, we set $\phi=\frac{|u|^{p(x)}}{f(v+\epsilon)}$. By the definition of solution \eqref{e32} we compute
\begin{align*}
\mu \int_\Omega a(x)f(v) \frac{|u|^{p(x)}}{f(v+\epsilon)} dx &\le \int_\Omega |\nabla_Xv|^{p-2}\nabla_X v \nabla_X\left(\frac{|u|^{p(x)}}{f(v+\epsilon)}\right) dx\\
& = \int_\Omega \left[|\nabla_X u|^{p(x)}-R(u,v+\epsilon)\right]dx\\
& = \int_\Omega|\nabla_X u|^{p(x)}dx -\int_\Omega L(u,v+\epsilon)dx.
\end{align*}
Taking the limit as $\epsilon \to 0^+$, applying Fatou's Lemma and Lebesgue dominated convergence theorem respectively on the left hand side and right hand side of the last expression, we obtain
\begin{align*}
0\le \int_\Omega |\nabla_X u|^{p(x)}- \mu \int_\Omega a(x)|u|^{p(x)}dx -\int_\Omega L(u,v)dx.
\end{align*}
Therefore we have 
\begin{align*}
0\le \int_\Omega|\nabla_X u|^{p(x)}dx- \mu \int_\Omega a(x)|u|^{p(x)} dx
\end{align*}
since $L(u,v)\ge 0$ almost everywhere in $\Omega$. This therefore completes the proof.
\qed

\begin{corollary}
Suppose there exists $\lambda>0$ and a strictly positive sup-solution of \eqref{e31} such that $\nabla_Xp(x)\nabla_Xv=0$. Then
\begin{align}\label{e36}
\int_\Omega|\nabla_X u|^{p(x)}dx \ge  \lambda \int_\Omega g(x)|u|^{p(x)} dx
\end{align}
for all $u\in W_0^{1,p(x)}(\Omega)$.
\end{corollary}
\proof
Applying Proposition \ref{Pro33} by setting $a(x)\equiv g(x)$, $\mu =\lambda$ and $f(v)=|v|^{p(x)-2}v$ , then one arrives at the conclusion \eqref{e35} at once.
\qed

\subsection{Principal frequency and domain monotonicity}
\begin{proposition}\label{Pro35}
Let there exists $\lambda$ and a strictly positive sup-solution $v \in W_0^{1,p(x)}(\Omega)$ of \eqref{e31} such that $\nabla_Xp(x)\nabla_Xv=0$.  Then we have
\begin{align}\label{e37}
\int_\Omega|\nabla_X u|^{p(x)}dx  \ge \lambda \int_\Omega g(x)|u|^{p(x)} dx
\end{align}
and
\begin{align}\label{e38}
\lambda_{1,p}(\Omega)\ge \lambda
\end{align}
for all $u \in W_0^{1,p(x)}(\Omega)$.
\end{proposition}

\proof
Suppose there exists $\lambda>0$, since $v$ is strictly positive sup-solution of \eqref{e31} in $\Omega$,  we have 
\begin{align}\label{e39}
\int_\Omega|\nabla_Xv|^{p(x)-2}\langle\nabla_Xv,\nabla_X\phi\rangle dx \ge \lambda \int_\Omega g(x)|v|^{p(x)-2} v\phi dx
\end{align} 
for all $\phi \in W_0^{1,p(x)}(\Omega)$.  For a given small $\epsilon>0$, setting $\disp \phi=\frac{|u|^{p(x)}}{(v+\epsilon)^{p(x)-1}}$ into \eqref{e39}. Then,  following the proof of the Proposition \ref{Pro33}, we arrive at \eqref{e37}.

Now, let $u_1 \in W_0^{1,p(x)}(\Omega)$ be the eigenfunction corresponding to the principal  eigenvalue $\lambda_{1,p}(\Omega)$.  We have
\begin{align}\label{e310}
\int_\Omega|\nabla_Xu_1|^{p(x)-2}\langle\nabla_Xu_1,\nabla_X\phi\rangle dx = \lambda_{1,p} \int_\Omega g(x)|u_1|^{p(x)-2} u_1\phi dx 
\end{align} 
for any $ \phi \in W_0^{1,p(x)}(\Omega)$.
Choosing $\epsilon>0$ (small) we can define via Picone identity that 
\begin{align}\label{e311}
0\le L(u_1,v+\epsilon)=R(u_1,v+\epsilon), \ \ v>0.
\end{align}
Integrating \eqref{e311} over $\Omega$ and then using \eqref{e39} with $\disp \phi=\frac{|u_1|^{p(x)}}{f(v+\epsilon)}$ and \eqref{e310} with $\phi=u_1$, we obtain
\begin{align*}
0& \le \int_\Omega L(u_1,v+\epsilon)dx = \int_\Omega R(u_1,v+\epsilon)dx\\
&=\int_\Omega|\nabla_X u_1|^{p(x)}dx - \int_\Omega \nabla_X\left(\frac{|u_1|^{p(x)}}{f(v+\epsilon)}\right) |\nabla_X v|^{p(x)-2}\nabla_X vdx\\
&=\int_\Omega|\nabla_X u_1|^{p(x)}dx + \int_\Omega  \frac{|u_1|^{p(x)}}{f(v+\epsilon)} \nabla_X^*( |\nabla_X v|^{p(x)-2}\nabla_X) vdx\\
&\le \lambda_{1,p}(\Omega)\int_\Omega g(x)|u_1|^{p(x)}dx - \lambda\int_\Omega g(x)\frac{|u_1|^{p(x)}}{f(v+\epsilon)}|v|^{p(x)-2}vdx.
\end{align*}
As usual, taking the limit as $\epsilon \to 0^+$, applying Fatou's Lemma and Lebesgue dominated convergence theorem, setting $f(v)=v^{p(x)-1}$,  we arrive at
\[0\le (\lambda_{1,p}(\Omega)-\lambda)\int_\Omega g(x)|u_1|^pdx,\]
which implies $\lambda_{1,p}(\Omega)\ge\lambda$.
\qed

As a corollary to the last proposition, we show strict monotonicity of the principal eigenvalue with respect to domain monotonicity. Let $\lambda_{1,p}(\Omega)>0$ be the principal eigenvalue of $\mathscr{L}_p$ on $\Omega$.
\begin{corollary}
 Suppose $\Omega_1\subset\Omega_2\subset\Omega$ and $\Omega_1\neq\Omega_2$. Let $u_1$ and $u_2$ be the eigenfunctions corresponding to $\lambda_{1,p}(\Omega_1)$ and $\lambda_{1,p}(\Omega_2)$ satisfying $\nabla_Xp(x)\nabla_Xu_1 =0$ and $\nabla_Xp(x)\nabla_Xu_2 =0$.  Then
\begin{align*}
\lambda_{1,p}(\Omega_1)> \lambda_{1,p}(\Omega_2)
\end{align*}
if they both exist.
\end{corollary}

\proof
Let $u_1$  and $u_2$ be positive eigenfunctions corresponding to $\lambda_{1,p}(\Omega_1)$ and $\lambda_{1,p}(\Omega_2)$, respectively. Clearly with $\phi\in C^\infty_0(\Omega)$, we have by Picone identity that
\[0\le \int_\Omega L(\phi,u_2)dx = \int_\Omega R(\phi,u_2)dx.\]
Replacing $\phi$ by $u_1$ and applying Proposition \ref{Pro35} we have 
\[\lambda_{1,p}(\Omega_1)- \lambda_{1,p}(\Omega_2)\ge 0.\]
If we have $\lambda_{1,p}(\Omega_1)= \lambda_{1,p}(\Omega_2)$, then  $L(u_1,u_2)=0$ a.e. in $\Omega$ and thus $u_1=\alpha u_2$ for some constant $\alpha>0$. However, this is   impossible when $\Omega_1\subset\Omega_2$ and $\Omega_1\neq\Omega_2$.
\qed

Next is the uniqueness and simplicity results.

\subsection{Uniqueness and simplicity of first eigenvalue}

\begin{proposition}\label{Pro37}
Let there exists $\lambda>0$ and a strictly positive solution $v \in W_0^{1,p(x)}(\Omega)$ of \eqref{e31} such that $\nabla_Xp(x)\nabla_Xv=0$. Then we have
\[\lambda_{1,p}(\Omega)= \lambda.\]
Moreover, let $u_1$ be the corresponding eigenfunction to $\lambda_{1,p}(\Omega)$. Then any other $u \in W_0^{1,p(x)}(\Omega)$ corresponding to $\lambda_{1,p}(\Omega)$ is a constant multiple of $u_1$.
\end{proposition}

\proof
Let $u_1\in W_0^{1,p(x)}(\Omega)$ be the eigenfunction corresponding to $\lambda_{1,p}(\Omega)$ and $u$ be a positive solution of \eqref{e31}.  Applying Picone identity by choosing $\epsilon>0$ (small) as follows:
\begin{align*}
0&\le \int_\Omega L(u,u_1+\epsilon) dx\\
&=\int_\Omega|\nabla_X u|^{p(x)}dx + \int_\Omega  \frac{u^{p(x)}}{f(u_1+\epsilon)} \nabla_X^*( |\nabla_X u_1|^{p(x)-2}\nabla_X) u_1dx\\
&= \lambda\int_\Omega g(x)|u|^{p(x)}dx -  \lambda_{1,p}(\Omega)\int_\Omega g(x)\frac{u^{p(x)}}{(u_1+\epsilon)^{p(x)-1}}|u_1|^{p(x)-2}u_1dx,
\end{align*}
where we have set $f(u_1+\epsilon)=(u_1+\epsilon)^{p(x)-1}$. Taking the limit as $\epsilon \to 0^+$, applying Fatou's Lemma and Lebesgue dominated convergence theorem,  then
\[\lambda_{1,p}(\Omega)\le\lambda.\]
On the other hand by Proposition \ref{Pro35}, we have 
\[\lambda_{1,p}(\Omega)\ge\lambda.\]
This therefore implies that $\lambda_{1,p}(\Omega)=\lambda$. By this we have proved the uniqueness part.

Now by the hypothesis of the theorem we have for $\phi,\psi \in C^\infty_0(\Omega)$ that
\begin{align}\label{e312}
\int_\Omega|\nabla_Xu|^{p(x)-2}\langle\nabla_Xu,\nabla_X\phi\rangle dx = \lambda_{1,p} \int_\Omega g(x)|u|^{p(x)-2} u\phi dx, 
\end{align}
\begin{align}\label{e313}
\int_\Omega|\nabla_Xu_1|^{p(x)-2}\langle\nabla_Xu_1,\nabla_X\psi\rangle dx = \lambda_{1,p} \int_\Omega g(x)|u_1|^{p(x)-2} u_1\psi dx.
\end{align}
Taking $\phi=u$ and $\psi=\frac{|u|^p}{(u_1+\epsilon)^{p-1}}$ into \eqref{e312} and \eqref{e313}, respectively, and sending $\epsilon \to 0^+$, we arrive at
\begin{align*}
\int_\Omega|\nabla_Xu|^{p(x)} dx& = \lambda_{1,p}\int_\Omega g(x)|u|^{p(x)} dx \\
&=\int_\Omega|\nabla_Xu_1|^{p(x)-2} \nabla_Xu_1  \nabla_X \Big(\frac{|u|^{p(x)}}{u_1^{p(x)-1}}\Big)dx,
\end{align*}
which implies (by choosing $f(u_1)=u_1^{p(x)-2}$)
\begin{align*}
\int_\Omega R(u,u_1)dx = \int_\Omega L(u,u_1)dx =0
\end{align*}
and consequently,  $\nabla_X(u/v)=0$, i.e.,  $u=\alpha u_1$  for some positive constant $\alpha>0$. 

\qed

The next proposition gives the sign changing nature of any other eigenfunction associated to an eigenvalue other than $\lambda_{1,p}(\Omega)$.
\begin{proposition}
Any eigenfunction $v$ corresponding to an eigenvalue $\lambda\neq \lambda_{1,p}(\Omega)$ such that $\nabla_Xp(x)\nabla_Xv=0$ changes sign.
\end{proposition}

\proof
By contradiction we suppose $v>0$ does not change sign (the case $v\le 0$ can be handled similarly). Let $\phi>0$ be an eigenfunction corresponding to $\lambda_{1,p}(\Omega)$. Choosing any $\epsilon>0$ as before,  applying Picone identity,  we have 
\begin{align*}
0&\le \int_\Omega L(\phi,v+\epsilon) dx\\
& =\int_\Omega\left[|\nabla_X \phi|^{p(x)} -\nabla_X\Big(\frac{\phi^{p(x)}}{f(v+\epsilon)}\Big)|\nabla_Xv|^{p(x)-2}\nabla_X v \right]dx\\
& =\int_\Omega |\nabla_X \phi|^{p(x)}dx + \int_\Omega \frac{\phi^{p(x)}}{f(v+\epsilon)} \mathscr{L}_pv dx.
\end{align*}
Since $\frac{\phi^{p(x)}}{(v+\epsilon)^{p(x)-1}}$ is admissible in the weak formulation of \eqref{e31} satisfied by $(\phi,\lambda)$, we arrive at 
\[0\le\lambda_{1,p}(\Omega)\int_\Omega g(x)|\phi|^{p(x)} dx - \lambda\int_\Omega \frac{\phi^{p(x)}}{f(v+\epsilon)} g(x) |v|^{p(x)-2}v dx.\]
Setting $f(v+\epsilon)=(v+\epsilon)^{p(x)-1}$ and letting $\epsilon \to 0^+$ in the last inequality as usual we obtain
\[0\le (\lambda_{1,p}-\lambda)\int_\Omega g(x) \phi^{p(x)} dx,\]
which is a contradiction since $\int_\Omega g(x)\phi^{p(x)}dx=1$. Thus $v$ must change sign.

\qed

\section{Variable exponent Caccioppoli estimates for general vector fields}
Picone identity is applied to prove some variable exponent Caccioppoli estimates for general vector fields in this section.  Recall that
$$1<p^-:=ess\inf_{x\in\bar{\Omega}} p(x)\le p(x) \le p^+:= ess\sup_{x\in\bar{\Omega}} p(x)<\infty.$$
Without giving rise to confusion but for simplicity sake we write $p:=p(x)$ and $q=:q(x)$.  We also denote $q^-:= ess\inf_{x\in\bar{\Omega}} q(x)$ and $q^+:= ess\sup_{x\in\bar{\Omega}} q(x)$.
\begin{theorem}\label{thm41}
Let $v$ be a positive sub-solution of \eqref{e31} in $\Omega\subset M$. Then for every fixed $q(x)>p(x)-1$, $p(x)>1$, $\nabla_Xv\nabla_Xp(x)=0$,  $\nabla_Xv\nabla_Xq(x)=0$  and $\lambda \in \mathbb{R}$, we have
\begin{align}
\int_\Omega v^{q-p}\phi^p|\nabla_Xv|^pdx \le C^{p^+}_{p,q}\int_\Omega v^q |\nabla_X\phi|^pdx + C_{\lambda,p,q} \int_\Omega g(x)v^q\phi^pdx
\end{align}
for every nonnegative functions $\phi\in C^\infty_0(\Omega)$, where\\
$$C^{p^+}_{p,q}:=\left( \frac{p^+}{q^--p^++1}\right)^{p^+} \ \ \ \text{and}\ \ \  C_{\lambda,p,q}:= 
\left( \frac{\lambda p^+}{q^--p^++1}\right).$$
\end{theorem}

\proof
Let $u=v^{q/p}\phi$, where $\phi$ is a nonnegative test function and $v$ is a sub-solution of \eqref{e31}, we compute
\begin{align*}
\nabla_X\left(v^{ q/p}\phi\right) &= \phi\nabla_X(v^{q/p})+v^{q/p}\nabla_X\phi\\
& = \phi v^{q/p}\ln v\left(\frac{\nabla_Xq}{p}-\frac{q\nabla_Xp}{p^2}\right)+\frac{q}{p}v^{\frac{q-p}{p}}\phi\nabla_Xv  + v^{q/p}\nabla_X\phi
\end{align*}
so that 
\begin{align*}
\langle\nabla_Xv, \nabla_X\left(v^{ q/p}\phi\right)\rangle & = \phi v^{q/p}\ln v\left(\frac{\nabla_Xq}{p}-\frac{q\nabla_Xp}{p^2}\right)\nabla_Xv \\
& \hspace{1cm} +\frac{q}{p}v^{\frac{q-p}{p}}\phi|\nabla_Xv|^2  + v^{q/p}\langle\nabla_X\phi,\nabla_Xv\rangle.
\end{align*}
Now using the the fact that $v$ is a sub-solution of \eqref{e31} and the condition that $\nabla_Xv\nabla_Xp(x)\equiv 0$ and $\nabla_Xv\nabla_Xq(x)\equiv 0$ in the Picone identity $L(u,v)\ge 0$, we have 
\begin{align}\label{e42}
0  & \le \int_\Omega L(v^{q/p}\phi,v)\nonumber\\
& = \int_\Omega |\nabla_X\left(v^{ q/p}\phi\right)|^pdx +\int_\Omega \frac{f'(v)}{(f(v))^2} |v^{q/p}|^p|\phi\nabla_Xv|^p dx\nonumber \\
& \ \ \ - \int_\Omega q \frac{|v^{q/p}\phi|^{p-1}}{f(v)}\phi v^{\frac{q-p}{p}} |\nabla_Xv|^pdx\\
& \ \ \ - \int_\Omega p \frac{|v^{q/p}\phi|^{p-1}}{f(v)} v^{q/p} |\nabla_Xv|^{p-2} \langle\nabla_X\phi,\nabla_Xv\rangle dx.\nonumber
\end{align}
Considering the condition $f'(v)\ge (p(x)-1)\left[f(v)^{\frac{p(x)-2}{p(x)-1}}\right]$, we can then choose $f(v)=v^{p(x)-1}$. Then \eqref{e42} reads
\begin{align}\label{e43}
0  \le &  \int_\Omega |\nabla_X\left(v^{ q/p}\phi\right)|^pdx +\int_\Omega (p-1)  v^{q-p} |\phi\nabla_Xv|^p dx - \int_\Omega q  v^{q-p} |\phi\nabla_Xv|^pdx \nonumber \\
& - \int_\Omega p |v^{\frac{q-p}{p}}\phi|^{p-1}  v^{q/p} |\nabla_Xv|^{p-2} \langle\nabla_X\phi,\nabla_Xv\rangle dx. 
\end{align}
Using the $\varepsilon(x)$-modified version of the Young's inequality in Lemma \ref{Y2} with $\Phi=v^{q/p}|\nabla_X\phi|$ and $\Psi=v^{\frac{q-p}{p}}\phi|\nabla_Xv|$, we can estimate the last term of \eqref{e43} as follows
\begin{align}\label{e44}
- \int_\Omega   p |v^{\frac{q-p}{p}}\phi|^{p-1}  & v^{q/p} |\nabla_Xv|^{p-2} \langle\nabla_X\phi,\nabla_Xv\rangle dx \nonumber \\
& \le \int_\Omega p |v^{\frac{q-p}{p}}\phi|^{p-1}   |\nabla_Xv|^{p-1} v^{q/p}  \nabla_X\phi dx \nonumber \\
&\le \int_\Omega \varepsilon^{1-p}v^q |\nabla_X\phi|^pdx +\int_\Omega \varepsilon(p-1) v^{q-p} |\phi\nabla_Xv|^pdx,
\end{align}
where $\varepsilon(x)$ is a continuous bounded function on $\Omega$, which will be chosen later. Substituting \eqref{e44} into \eqref{e43} we get 
\begin{align*}
0  & \le \int_\Omega |\nabla_X\left(v^{ q/p}\phi\right)|^pdx - \int_\Omega [q-p+1-\varepsilon(p-1)] v^{q-p} |\phi\nabla_Xv|^pdx\\
& \hspace{1cm} + \int_\Omega \varepsilon^{1-p}v^q |\nabla_X\phi|^pdx \\
& \le \lambda  \int_\Omega g(x) |v^{ q/p}\phi|^pdx  - \mathcal{C}^1_{\epsilon,p,q}  \int_\Omega  v^{q-p} |\phi\nabla_Xv|^pdx + \mathcal{C}^2_{\epsilon,p}  \int_\Omega  v^q |\nabla_X\phi|^pdx, 
\end{align*}
where we have used  $\disp \int_\Omega|\nabla_Xu|^{p(x)}dx \le \lambda \int_\Omega g(x)|u|^{p(x)}dx$ for the sub-solution of \eqref{e31}. Here 
\begin{align*}
\mathcal{C}^1_{\epsilon,p,q} := q^--p^++1-\bar{\varepsilon}(p^+-1)
\  \ \text{and}\ \ \mathcal{C}^2_{\epsilon,p} : = \bar{\varepsilon}^{1-p^+},
\end{align*}
where $\bar{\varepsilon}:=\sup_\Omega\varepsilon(x)$.  

\noindent Rearranging the last inequality we arrive at 
\begin{align*}
 \int_\Omega  v^{q-p} |\phi\nabla_Xv|^pdx \le  \frac{ \mathcal{C}^2_{\epsilon,p}}{\mathcal{C}^1_{\epsilon,p,q} }  \int_\Omega  v^q |\nabla_X\phi|^pdx + \frac{\lambda}{\mathcal{C}^1_{\epsilon,p,q} }  \int_\Omega g(x) |v^{ q/p}\phi|^pdx. 
\end{align*}
We can now choose a suitable number $\bar{\varepsilon}$ as  $\disp \bar{\varepsilon}:= \frac{q^--p^++1}{p^+}$ and then compute
\begin{align*}
 \frac{1}{\mathcal{C}^1_{\epsilon,p,q} } & :=  \frac{1}{q^--p^++1-\bar{\varepsilon}(p^+-1)}   =  \frac{p^+}{q^--p^++1},
 \\
  \frac{ \mathcal{C}^2_{\epsilon,p}}{\mathcal{C}^1_{\epsilon,p,q} } &: = \frac{ \bar{\varepsilon}^{1-p^+}}{q^--p^++1-\bar{\varepsilon}(p^+-1)} =  \left( \frac{p^+}{q^--p^++1}\right)^{p^+}.
\end{align*}
The proof is therefore complete.
\qed

The following two corollaries can be deduced from Theorem \ref{thm41} using the same assumptions.
\begin{corollary}\label{cor42}
Let $v$ be a positive sub-solution of \eqref{e31} in $\Omega$ satisfying $\nabla_Xv\nabla_Xp(x)=0$. If $g(x)\equiv 0$ and $p(x)=q(x)$ in $\Omega$. Then we have 
\begin{align*}
\int_\Omega \phi^{p(x)}|\nabla_Xv|^{p(x)}dx  \le (p^+)^{p^+} \int_\Omega v^{p(x)} |\nabla_X\phi|^{p(x)} dx
\end{align*}
for every nonnegative function $\phi \in C^\infty_0(\Omega)$.
\end{corollary}

\begin{corollary}\label{cor43}
Let $v$ be a positive sub-solution of \eqref{e31} in $\Omega$ satisfying $\nabla_Xv\nabla_Xp(x)=0$.  Letting $\lambda =1$ and $p(x)=q(x)$ in $\Omega$. Then we have 
\begin{align*}
\int_\Omega \phi^{p(x)}|\nabla_Xv|^{p(x)}dx  \le (p^+)^{p^+} \int_\Omega v^{p(x)} |\nabla_X\phi|^{p(x)} dx + p^+\int_\Omega g(x)  v^{q(x)} \phi^{p(x)} dx
\end{align*}
for every nonnegative function $\phi \in C^\infty_0(\Omega)$.
\end{corollary}

\begin{remark}
Suppose $M=\mathbb{R}^n$, $p(x)= p$ (constant) and $q(x)=q$ (constant):
\begin{enumerate}[(i)]
\item Corollary \ref{cor42} reduces to \cite[Corollary 3.1]{Ja4} and \cite[equation 5.27]{PRS}.
\item Corollary \ref{cor42} reduces to \cite[Corollary A.6]{LLM1}.
\end{enumerate}
\end{remark}

We remark also that analogous result to Theorem \ref{thm41} holds for positive sup-solutions of \eqref{e31}  with $q(x) < p(x)-1$. 
\begin{theorem}\label{thm43}
Let $v$ be a positive sup-solution of \eqref{e31} in $\Omega\subset M$.  Then for every fixed $q(x)<p(x)-1$, $p(x)>1$, $\nabla_Xv\nabla_Xp(x)=0$,  $\nabla_Xv\nabla_Xq(x)=0$  and $\lambda \in \mathbb{R}$, we have
\begin{align}\label{a2}
\int_\Omega v^{q-p}\phi^p|\nabla_Xv|^pdx \le C^{p^-}_{p,q}\int_\Omega v^q |\nabla_X\phi|^pdx + C_{\lambda,p,q} \int_\Omega g(x)v^q\phi^pdx
\end{align}
for every nonnegative functions $\phi\in C^\infty_0(\Omega)$, where\\
$$C^{p^-}_{p,q}:=\left( \frac{p^+}{p^--q^+-1}\right)^{p^+} \ \ \ \text{and}\ \ \  C^-_{\lambda,p,q}:= -
\left( \frac{\lambda p^+}{p^--q^+-1}\right).$$
\end{theorem}

\begin{remark}
Setting $q=0$  in \eqref{a2} we obtain a particular case whose right hand side is independent of the nonnegative function $v$. That is
\begin{align}\label{a3}
\int_\Omega |\phi\nabla_X\log v|^pdx \le \left(\frac{p^+}{p^--1}\right)^{p^+}\int_\Omega  |\nabla_X\phi|^pdx - \left( \frac{\lambda p^+}{p^--1}\right) \int_\Omega g(x)\phi^pdx.
\end{align}
This is the variable exponent logarithmic Caccioppolli inequality.  Precisely,  If $g(x)\equiv 0$, then \eqref{a3} reduces to a new version of the well known  logarithmic Caccioppolli inequality for positive $p(x)$-superharmonic functions
\begin{align*}
\int_\Omega |\phi\nabla_X\log v|^{p(x)}dx \le \left(\frac{p^+}{p^--1}\right)^{p^+}\int_\Omega  |\nabla_X\phi|^{p(x)}dx,
\end{align*}
where $1<p^-<p^+<\infty$.  Note that $v \in W^{1,p(x)}_{\text{loc}}$ is
said to be $p(x)$-superharmonic if it satisfies $\int_\Omega|\nabla_Xu|^{p(x)-2}\langle\nabla_Xu,\nabla_X\phi\rangle dx \ge 0$. Interested reader is hereby referred to \cite{LindqM}  and \cite{Lindq} for $p$($=$ constant)-superharmonic case. 
\end{remark}

\section*{Compliance with Ethical Standards}

\subsection*{Conflicts of Interest}  The authors declare that they have no conflict of interests.
 
\subsection*{Ethical Approval } This article does not contain any studies with human participants or animals performed by any of the authors.

\subsection*{Funding} This project does not receive any funding.

\subsection*{Acknowledgement}
This paper was completed during the first author's research visit to Ghent Analysis and PDE Centre, Ghent University.  He therefore gratefully acknowledges the  research supports of  IMU-Simons African Fellowship Grant  and EMS-Simons for African  program.  He also thanks his host Professor Michael Ruzhansky for useful discussions on this project.  Part of the results of this paper have been presented in \cite{A2} (2022 Proceedings of International E-Conference on Mathematical and Statistical Sciences: https://icomss22.selcuk.edu.tr).


\begin{thebibliography}{20}

\bibitem{A1}
 		A. Abolarinwa,  
 		Anisotropic Picone type identities for general vector fields and some Applications.  In: Ruzhansky, M., Van Bockstal, K. (eds) Extended Abstracts 2021/2022. APDEGS 2021. Trends in Mathematics, vol 2. Birkhäuser, Cham (2024).


\bibitem{A2}
		A. Abolarinwa, 
		{\it Generalised Picone identity and first eigenvalue for $p(x)$-sub-Laplacian on stratified groups},   2022 Proceedings of International E-Conference on Mathematical and Statistical Sciences: https://icomss22.selcuk.edu.tr


\bibitem{AMS}
		E. Acerbi, G. Mingione, G. A. Seregin,
		{\it Regularity results for parabolic systems related to a class of non-Newtonian fluids},
 		Ann. Inst. H. Poincar\'e Anal.  Non  Line\'aire, 21(1) (2004),  25--60.

\bibitem{AMSo}
		R.  Aboulaich, D. Meskine, A. Souissi,
		{\it New diffusion models in image processing}, 
		Comput.  Math. Appl., 56(4) (2008),  874--882.
		
\bibitem{Al1} 
		W. Allegretto, 
		{\it Positive solutions and spectral properties of weakly coupled elliptic systems}, 
		J.  Math.  Anal.  Appl., 120 (1986), 723--729.
		
\bibitem{AH} 
		W. Allegretto, Y. X. Huang, 
		{\it A Picone's identity for the $p$-Laplacian and applications}, 
		Nonl. Anal. Th. Meth. Appl., 32(7) (1998),  819--830.


\bibitem{Al3} 
		W. Allegretto, 
		{\it Form estimates for the $p(x)$-Laplacian}, 
		Proc.  Amer.  Math.  Soc., 135(7) (2007),  2177--2185.
		
\bibitem{Alv}
		C. O. Alves,
		{\it Existence of solution for a degenerate $p(x)$-Laplacian equation in $\mathbb{R}^N$},
		 J.  Math.  Anal. Appl. , 345(2) (2008),  731--742.

		 
\bibitem{Ba}
		K. Bal, 
		{\it Generalized Picone's identity and its applications}, 
		Electron.  J. Diff. Equ.,  2013(243), (2013), pp. 1--6. 


\bibitem{CFRW}
		D. Cruiz-Uribe, A. Fiorenza, M. Ruzhansky, J. Wirth,
		{\it Variable Lebesgue spaces and Hyperbolic systems},
 		Springer Basel, 2014.
		
				
\bibitem{Deng}
		S-G. Deng,
		{\it Eigenvalues of $p(x)$-Laplacian Dirichlet Steklov problem}, 
		J.  Math.  Anal. Appl., 339(2) (2008),  925--937.

\bibitem{DHHR}
		L. Diening, P.  Harjulehto,  P. H\"ast\"o, \'Ut V. L\^e,  M.   Ruzicka,
		{\it Lebesgue and Sobolev spaces with variable exponents},
		 Springer 2011.	
		 
		 		
\bibitem{DT}
 			G. Dwivedi,  J. Tyagi, 
 			{\it Picone's identity for biharmonic operators on Heisenberg group and its applications}, 				
 			Nonlinear Differ. Equ. Appl.,  23(2) (2016),  1--26.		

\bibitem{Fan}
		X. Fan, 
		{\it Eigenvalues of the $p(x)$-Laplacian Neumann problems}, 
		Nonlinear Anal., 67 (2007),  2982--2992.
		
\bibitem{FZ}
		X. Fan, D. Zhao, 
		{\it On the spaces $L^{p(x)}(\Omega)$ and $W^{m,p(x)}(\Omega)$}, 
		J.  Math.  Anal. Appl.,  263(2) (2001),  424--446.
		
\bibitem{FZZ}
		X. Fan, Q. Zhang, D. Zhao,
		{\it Eigenvalues of $p(x)$-Laplacian Dirichlet problem}, 
		J.  Math.  Anal. Appl.,  302(2) (2005),  306--317.

				
\bibitem{Fe}
		T. Feng,
		{\it A new nonlinear Picone identity and applications}, 
		Mathematica Applicata, 30(2) (2017),  278--283.
		
\bibitem{FH}
		T. Feng, J. Han,
		{\it A new variable exponent Picone identity and applications},
		Math. Ineq. Appl., 22(1) (2019), 65--75.
		

\bibitem{FL}
		G. Franzina, P.  Lindqvist,
		{\it An eigenvalue problem with variable exponents},
		Nonlinear Anal.,  85 (2013), 1--16.
		
\bibitem{HHV}
		P.  Harjulehto,  P. H\"ast\"o, \'Ut V. L\^e,  M.   Nuortio,
		{\it Overview of differential equations with non-standard growth},
		 Nonlinear Anal.,  72 (2010),  4551--4574.		
		
\bibitem{Hom}
		H.  H\"ormander, 
		{\it Hypoelliptic second-order differential equations},
		 Acta Math.,  43 (1967),  147--171.
		 	

\bibitem{Ja3}
		J. Jaro\v{s},
		{\it A-harmonic Picone's identity with applications},
		 Annali di Matematica, 194(3) (2015), 719--729.

\bibitem{Ja4}
		J. Jaro\v{s},
		{\it Caccioppoli estimates through an anisotropic Picone's identity},
		Proc. Amer.  Math. Soc., 143 (2015), 1137--1144.
	
	
\bibitem{LindqM}
		P. Lindqvist, J. J. Manfredi,
		{\it Hardy's inequality from a logarithmic Caccioppoli Estimate}, 
		arxiv:
\bibitem{Lindq}
		P. Lindqvist, 
		{\it On the definition and properties of $p$-superharmonic functions}, 
		Journal für die reine und angewandte Mathematik, 
		1986 (365), 1986,  67 -- 79. 
			
\bibitem{LLM1}
		V. Liskevich, S. Lyakhova,  V. Moroz, 
		{\it Positive solutions to nonlinear $p$-Laplace equations with Hardy potential in exterior domains}, 
		J.  Diff. Eq.  232(1) (2007),  212 -- 252.
	
	
\bibitem{MV}
		M.  Mihailescu, C. Varga,
		{\it Multiplicity results for some elliptic problems with nonlinear boundary conditions 					involving variable exponents},  
		Comput.  Math.  Appl., 62(9) (2011),  3464--3471.
		
\bibitem{NZW}
		P. Niu, H. Zhang, Y. Wang,
		{\it Hardy Type and Rellich Type Inequalities on the Heisenberg Group},
		Proc. Amer. Math. Soc., 129(12) (2001), 3623--3630.
		
\bibitem{PRS}
		S. Pigola,  M. Rigoli,  A. G. Setti, 
		{\it Vanishing and finiteness results in geometric analysis}, Progress in Mathematics, vol. 266, Birkh\"auser Verlag, Basel, 2008. A generalization of the Bochner technique
										
\bibitem{RS}
		  M. Ruzhansky, D. Suragan,
		{\it Hardy inequalities on homogeneous groups}, 
		Progress in Math. Vol. 327, 	Birkh\"auser, 588 pp, (2019).
		
		
\bibitem{RSS1}
		M. Ruzhansky, B. Sabitbek, D. Suragan,
		{\it Principal frequency of $p$-versions of sub-Laplacians for general 										vector fields},  
 		Z. Anal. Anwend.,   40 (2021), 97--109.		
		
\bibitem{RSS2}
		M. Ruzhansky, B. Sabitbek, D. Suragan,
		{\it Weighted anisotropic Hardy and Rellich type inequalities for general 				vector fields},  
		Nonlinear Differ. Eqn. Appl.,  26, 13 (2019). 
		
\bibitem{RSS3}
 			M. Ruzhansky, B. Sabitbek, D. Suragan,
 		{\it Weighted $L^p$-Hardy and $L^p$-Rellich inequalities with boundary terms on stratified Lie groups},
 		Rev. Mat. Complut.,  32(1) (2019),  19--35.

\bibitem{Ru}
		M.  Ruzicka,
		{\it Electrorheological Fluids Modeling and Mathematical Theory}, 
		Springer-Verlag,  Berlin,  2000.

		 
 \bibitem{SY}	 
		D. Surugan, N. Yessirkegenov,
		{\it Generalised nonlinear Picone identities for $p$-sub-Laplacians and $p$-biharmonic operators and applications},
		Adv. Op. Th., 6(53) (2021), 1--17.		

\bibitem{Yo1}
		N. Yoshida,
		{\it Picone identity for quasilinear elliptic equations with $p(x)$-Laplacians and 							Sturmianian comparison theory},
		Appl.  Math.  Comput.,  225 (1) (2013),  79--91.
	
\bibitem{Yo2}
		N. Yoshida,
		{\it Picone-type inequality and Sturmian comparison theorems for quasilinear elliptic 				operators with $p(x)$-Laplacians},
		Electron.  J.  Differ.  Equ.,  2012(01), (2012), 1--9.		

	
\bibitem{Yo3}
		N. Yoshida,
		{\it Picone identities for half-linear elliptic operators with $p(x)$-Laplacians and 						applications to Sturmian comparison theory},
		 Nonlinear Anal.,  74 (2011),  5631--5642.

\bibitem{Tir}
		A. Tirayaki,
		{\it Generalized nonlinear Picone's identity for the $p$-Laplacian and its 									applications},
		Electron.  J.  Diff. Equ.,  2016(269),  (2016), 1--7.
				
\bibitem{Ty}
		J. Tyagi, 
		{\it A nonlinear picone's identity and its applications}, 
		Appl.  Math.  Letters,  26 (2013), 624--626.				
\end{thebibliography}
\end{document}